\documentclass[12pt]{amsart}
\usepackage{amssymb}
\begin{document}
\theoremstyle{plain}
\newtheorem{thm}{Theorem}[section]
\newtheorem{theorem}[thm]{Theorem}
\newtheorem{lemma}[thm]{Lemma}
\newtheorem{corollary}[thm]{Corollary}
\newtheorem{proposition}[thm]{Proposition}
\newtheorem{addendum}[thm]{Addendum}
\newtheorem{variant}[thm]{Variant}
\theoremstyle{definition}
\newtheorem{notations}[thm]{Notations}
\newtheorem{question}[thm]{Question}
\newtheorem{problem}[thm]{Problem}
\newtheorem{remark}[thm]{Remark}
\newtheorem{remarks}[thm]{Remarks}
\newtheorem{definition}[thm]{Definition}
\newtheorem{claim}[thm]{Claim}
\newtheorem{assumption}[thm]{Assumption}
\newtheorem{assumptions}[thm]{Assumptions}
\newtheorem{properties}[thm]{Properties}
\newtheorem{example}[thm]{Example}
\numberwithin{equation}{thm}
\catcode`\@=11
\def\opn#1#2{\def#1{\mathop{\kern0pt\fam0#2}\nolimits}}
\def\bold#1{{\bf #1}}%
\def\underrightarrow{\mathpalette\underrightarrow@}
\def\underrightarrow@#1#2{\vtop{\ialign{$##$\cr
 \hfil#1#2\hfil\cr\noalign{\nointerlineskip}%
 #1{-}\mkern-6mu\cleaders\hbox{$#1\mkern-2mu{-}\mkern-2mu$}\hfill
 \mkern-6mu{\to}\cr}}}
\let\underarrow\underrightarrow
\def\underleftarrow{\mathpalette\underleftarrow@}
\def\underleftarrow@#1#2{\vtop{\ialign{$##$\cr
 \hfil#1#2\hfil\cr\noalign{\nointerlineskip}#1{\leftarrow}\mkern-6mu
 \cleaders\hbox{$#1\mkern-2mu{-}\mkern-2mu$}\hfill
 \mkern-6mu{-}\cr}}}
\let\amp@rs@nd@\relax
\newdimen\ex@
\ex@.2326ex
\newdimen\bigaw@
\newdimen\minaw@
\minaw@16.08739\ex@
\newdimen\minCDaw@
\minCDaw@2.5pc
\newif\ifCD@
\def\minCDarrowwidth#1{\minCDaw@#1}
\newenvironment{CD}{\@CD}{\@endCD}
\def\@CD{\def\A##1A##2A{\llap{$\vcenter{\hbox
 {$\scriptstyle##1$}}$}\Big\uparrow\rlap{$\vcenter{\hbox{%
$\scriptstyle##2$}}$}&&}%
\def\V##1V##2V{\llap{$\vcenter{\hbox
 {$\scriptstyle##1$}}$}\Big\downarrow\rlap{$\vcenter{\hbox{%
$\scriptstyle##2$}}$}&&}%
\def\={&\hskip.5em\mathrel
 {\vbox{\hrule width\minCDaw@\vskip3\ex@\hrule width
 \minCDaw@}}\hskip.5em&}%
\def\verteq{\Big\Vert&&}%
\def\noarr{&&}%
\def\vspace##1{\noalign{\vskip##1\relax}}\relax\let\amp@rs@nd@&\iffalse}\fi
 \CD@true\vcenter\bgroup\relax\let\\=\cr\iffalse}\fi\tabskip\z@skip\baselineskip20\ex@
 \lineskip3\ex@\lineskiplimit3\ex@\halign\bgroup
 &\hfill$\m@th##$\hfill\cr}
\def\@endCD{\cr\egroup\egroup}
\def\>#1>#2>{\amp@rs@nd@\setbox\z@\hbox{$\scriptstyle
 \;{#1}\;\;$}\setbox\@ne\hbox{$\scriptstyle\;{#2}\;\;$}\setbox\tw@
 \hbox{$#2$}\ifCD@
 \global\bigaw@\minCDaw@\else\global\bigaw@\minaw@\fi
 \ifdim\wd\z@>\bigaw@\global\bigaw@\wd\z@\fi
 \ifdim\wd\@ne>\bigaw@\global\bigaw@\wd\@ne\fi
 \ifCD@\hskip.5em\fi
 \ifdim\wd\tw@>\z@
 \mathrel{\mathop{\hbox to\bigaw@{\rightarrowfill}}\limits^{#1}_{#2}}\else
 \mathrel{\mathop{\hbox to\bigaw@{\rightarrowfill}}\limits^{#1}}\fi
 \ifCD@\hskip.5em\fi\amp@rs@nd@}
\def\<#1<#2<{\amp@rs@nd@\setbox\z@\hbox{$\scriptstyle
 \;\;{#1}\;$}\setbox\@ne\hbox{$\scriptstyle\;\;{#2}\;$}\setbox\tw@
 \hbox{$#2$}\ifCD@
 \global\bigaw@\minCDaw@\else\global\bigaw@\minaw@\fi
 \ifdim\wd\z@>\bigaw@\global\bigaw@\wd\z@\fi
 \ifdim\wd\@ne>\bigaw@\global\bigaw@\wd\@ne\fi
 \ifCD@\hskip.5em\fi
 \ifdim\wd\tw@>\z@
 \mathrel{\mathop{\hbox to\bigaw@{\leftarrowfill}}\limits^{#1}_{#2}}\else
 \mathrel{\mathop{\hbox to\bigaw@{\leftarrowfill}}\limits^{#1}}\fi
 \ifCD@\hskip.5em\fi\amp@rs@nd@}
\newenvironment{CDS}{\@CDS}{\@endCDS}
\def\@CDS{\def\A##1A##2A{\llap{$\vcenter{\hbox
 {$\scriptstyle##1$}}$}\Big\uparrow\rlap{$\vcenter{\hbox{%
$\scriptstyle##2$}}$}&}%
\def\V##1V##2V{\llap{$\vcenter{\hbox
 {$\scriptstyle##1$}}$}\Big\downarrow\rlap{$\vcenter{\hbox{%
$\scriptstyle##2$}}$}&}%
\def\={&\hskip.5em\mathrel
 {\vbox{\hrule width\minCDaw@\vskip3\ex@\hrule width
 \minCDaw@}}\hskip.5em&}
\def\verteq{\Big\Vert&}
\def\novarr{&}
\def\noharr{&&}
\def\SE##1E##2E{\slantedarrow(0,18)(4,-3){##1}{##2}&}
\def\SW##1W##2W{\slantedarrow(24,18)(-4,-3){##1}{##2}&}
\def\NE##1E##2E{\slantedarrow(0,0)(4,3){##1}{##2}&}
\def\NW##1W##2W{\slantedarrow(24,0)(-4,3){##1}{##2}&}
\def\slantedarrow(##1)(##2)##3##4{%
\thinlines\unitlength1pt\lower 6.5pt\hbox{\begin{picture}(24,18)%
\put(##1){\vector(##2){24}}%
\put(0,8){$\scriptstyle##3$}%
\put(20,8){$\scriptstyle##4$}%
\end{picture}}}
\def\vspace##1{\noalign{\vskip##1\relax}}\relax\let\amp@rs@nd@&\iffalse}\fi
 \CD@true\vcenter\bgroup\relax\let\\=\cr\iffalse}\fi\tabskip\z@skip\baselineskip20\ex@
 \lineskip3\ex@\lineskiplimit3\ex@\halign\bgroup
 &\hfill$\m@th##$\hfill\cr}
\def\@endCDS{\cr\egroup\egroup}
\newdimen\TriCDarrw@
\newif\ifTriV@
\newenvironment{TriCDV}{\@TriCDV}{\@endTriCD}
\newenvironment{TriCDA}{\@TriCDA}{\@endTriCD}
\def\@TriCDV{\TriV@true\def\TriCDpos@{6}\@TriCD}
\def\@TriCDA{\TriV@false\def\TriCDpos@{10}\@TriCD}
\def\@TriCD#1#2#3#4#5#6{%
\setbox0\hbox{$\ifTriV@#6\else#1\fi$}
\TriCDarrw@=\wd0 \advance\TriCDarrw@ 24pt
\advance\TriCDarrw@ -1em
\def\SE##1E##2E{\slantedarrow(0,18)(2,-3){##1}{##2}&}
\def\SW##1W##2W{\slantedarrow(12,18)(-2,-3){##1}{##2}&}
\def\NE##1E##2E{\slantedarrow(0,0)(2,3){##1}{##2}&}
\def\NW##1W##2W{\slantedarrow(12,0)(-2,3){##1}{##2}&}
\def\slantedarrow(##1)(##2)##3##4{\thinlines\unitlength1pt
\lower 6.5pt\hbox{\begin{picture}(12,18)%
\put(##1){\vector(##2){12}}%
\put(-4,\TriCDpos@){$\scriptstyle##3$}%
\put(12,\TriCDpos@){$\scriptstyle##4$}%
\end{picture}}}
\def\={\mathrel {\vbox{\hrule
   width\TriCDarrw@\vskip3\ex@\hrule width
   \TriCDarrw@}}}
\def\>##1>>{\setbox\z@\hbox{$\scriptstyle
 \;{##1}\;\;$}\global\bigaw@\TriCDarrw@
 \ifdim\wd\z@>\bigaw@\global\bigaw@\wd\z@\fi
 \hskip.5em
 \mathrel{\mathop{\hbox to \TriCDarrw@
{\rightarrowfill}}\limits^{##1}}
 \hskip.5em}
\def\<##1<<{\setbox\z@\hbox{$\scriptstyle
 \;{##1}\;\;$}\global\bigaw@\TriCDarrw@
 \ifdim\wd\z@>\bigaw@\global\bigaw@\wd\z@\fi
 \mathrel{\mathop{\hbox to\bigaw@{\leftarrowfill}}\limits^{##1}}
 }
 \CD@true\vcenter\bgroup\relax\let\\=\cr\iffalse}\fi
 \tabskip\z@skip\baselineskip20\ex@
 \lineskip3\ex@\lineskiplimit3\ex@
 \ifTriV@
 \halign\bgroup
 &\hfill$\m@th##$\hfill\cr
#1&\multispan3\hfill$#2$\hfill&#3\\
&#4&#5\\
&&#6\cr\egroup%
\else
 \halign\bgroup
 &\hfill$\m@th##$\hfill\cr
&&#1\\%
&#2&#3\\
#4&\multispan3\hfill$#5$\hfill&#6\cr\egroup
\fi}
\def\@endTriCD{\egroup}
\newcommand{\sA}{{\mathcal A}}
\newcommand{\sB}{{\mathcal B}}
\newcommand{\sC}{{\mathcal C}}
\newcommand{\sD}{{\mathcal D}}
\newcommand{\sE}{{\mathcal E}}
\newcommand{\sF}{{\mathcal F}}
\newcommand{\sG}{{\mathcal G}}
\newcommand{\sH}{{\mathcal H}}
\newcommand{\sI}{{\mathcal I}}
\newcommand{\sJ}{{\mathcal J}}
\newcommand{\sK}{{\mathcal K}}
\newcommand{\sL}{{\mathcal L}}
\newcommand{\sM}{{\mathcal M}}
\newcommand{\sN}{{\mathcal N}}
\newcommand{\sO}{{\mathcal O}}
\newcommand{\sP}{{\mathcal P}}
\newcommand{\sQ}{{\mathcal Q}}
\newcommand{\sR}{{\mathcal R}}
\newcommand{\sS}{{\mathcal S}}
\newcommand{\sT}{{\mathcal T}}
\newcommand{\sU}{{\mathcal U}}
\newcommand{\sV}{{\mathcal V}}
\newcommand{\sW}{{\mathcal W}}
\newcommand{\sX}{{\mathcal X}}
\newcommand{\sY}{{\mathcal Y}}
\newcommand{\sZ}{{\mathcal Z}}
\newcommand{\A}{{\mathbb A}}
\newcommand{\B}{{\mathbb B}}
\newcommand{\C}{{\mathbb C}}
\newcommand{\D}{{\mathbb D}}
\newcommand{\E}{{\mathbb E}}
\newcommand{\F}{{\mathbb F}}
\newcommand{\G}{{\mathbb G}}
\newcommand{\HH}{{\mathbb H}}
\newcommand{\I}{{\mathbb I}}
\newcommand{\J}{{\mathbb J}}
\renewcommand{\L}{{\mathbb L}}
\newcommand{\M}{{\mathbb M}}
\newcommand{\N}{{\mathbb N}}
\renewcommand{\P}{{\mathbb P}}
\newcommand{\Q}{{\mathbb Q}}
\newcommand{\R}{{\mathbb R}}
\newcommand{\T}{{\mathbb T}}
\newcommand{\U}{{\mathbb U}}
\newcommand{\V}{{\mathbb V}}
\newcommand{\W}{{\mathbb W}}
\newcommand{\X}{{\mathbb X}}
\newcommand{\Y}{{\mathbb Y}}
\newcommand{\Z}{{\mathbb Z}}
\newcommand{\id}{{\rm id}}
\newcommand{\rank}{{\rm rank}}
\newcommand{\sEnd}{{\mathcal End}}
\newcommand{\End}{{\rm End}}
\title[Families of abelian varieties]{Families of
abelian varieties over curves with maximal Higgs field}
\author[Eckart Viehweg]{Eckart Viehweg}
\address{Universit\"at Essen, FB6 Mathematik, 45117 Essen, Germany}
\email{ viehweg@uni-essen.de}
\thanks{This work has been supported by the ``DFG-Schwerpunktprogramm
Globale Methoden in der Komplexen Geometrie''. The second named author
is supported by a grant from the Research
Grants Council of the Hong Kong
Special Administrative Region, China
(Project No. CUHK 4239/01P)}
\author[Kang Zuo]{Kang Zuo}
\address{The Chinese University of Hong Kong, Department of Mathematics,
Shatin, Hong Kong}
\email{kzuo@math.cuhk.edu.hk}
\maketitle
Throughout this note, $Y$ will denote a non-singular complex projective
curve, and $f:A\to Y$ a family of abelian varieties, with
$A$ non-singular.
We write $U\subset Y$ for an open dense subscheme, with
$$
f:A_0=f^{-1}(U)\>>> U
$$
smooth, $S = Y\setminus U$, and $\Delta=f^{-1}(S)$. Consider the
weight $1$ variation of Hodge structures given by $f:A_0\to U$,
i.e. $R^1f_*\Z_{A_0}$. We will assume throughout this note,
that the monodromy of $R^1f_*\Z_{A_0}$ around all points in $S$ is
unipotent. We write
$$
(E, \theta)=(E^{1,0}\oplus E^{0,1}, \theta_{1,0})
$$
for the Higgs-bundles induced by the Deligne extension of
$(R^1f_*\Z_{A_0})\otimes \sO_U$. Hence
$E^{1,0}=f_*\Omega^1_{A/Y}(\log \Delta)$ and
$E^{0,1}=R^1f_*\sO_A$.
The Higgs field is given by the edge morphisms
$$
f_*\Omega^1_{A/Y}(\log \Delta)\to R^1f_*\sO_A \otimes
\Omega^{1}_{Y}(\log S)
$$
of the tautological sequence
$$
0\to {f}^*\Omega^1_Y(\log S)\to \Omega^1_{A}(\log \Delta) \to
\Omega_{A/Y}^1(\log \Delta))\to 0.
$$
By \cite{Kol} $E^{1,0}$ is a direct sum $F^{1,0} \oplus N^{1,0}$ with
$F^{1,0}$ ample and $N^{1,0}={\rm Ker}(\theta_{1,0})$ flat.

Correspondingly, we have $E^{0,1}=F^{0,1}\oplus N^{0,1}$ and $E$
is the direct sum of the Higgs bundles
\begin{equation}\label{splitt1}
(F=F^{1,0}\oplus F^{0,1}, \theta_{1,0})  \mbox{ \ \ \ and \ \ \ }
(N^{1,0}\oplus N^{0,1},0).
\end{equation}

For $g_0={\rm rank}(F^{1,0})$ the Arakelov inequalities
(\cite{Pet}, \cite{J-Z}) say
\begin{equation}\label{arakelov_ineq}
2\cdot \deg(F^{1,0}) \leq g_0\cdot \deg(\Omega^1_Y(\log S)).
\end{equation}
In this note we will try to understand the geometry of families
$f:A\to Y$, for which (\ref{arakelov_ineq}) is an equality, or as
we will say shortly, of families reaching the Arakelov bound.

A family of abelian varieties is reaching the Arakelov bound
if and only if the Higgs field is maximal (see \ref{semistable}, i), i.e. if
$\theta_{1,0}:F^{1,0} \to F^{0,1}\otimes \Omega^1_Y(\log S)$
is an isomorphism.

For families of elliptic curves, the maximality of the Higgs field
implies that the family is modular:

\begin{proposition}\label{modular}
Let $h:E \to Y$ be a semi stable family
of elliptic curves, smooth over $U\subset Y$ with $U\neq Y$.
If $E\to Y$ is reaching the Arakelov bound, $E\to Y$ is modular,
i.e. $U$ is the quotient of the upper half plane $\sH$ by a subgroup
of $\text{Sl}_2(\Z)$ of finite index, and the morphism $U \to \C=\sH/
\text{Sl}_2(\Z)$ is given by the $j$-invariant of the fibres.
\end{proposition}

In the higher dimensional case one could hope, that the general fibre
of a family with maximal Higgs field is quite special, and that the base
curve is again a Shimura curve. The corresponding question has been considered
in \cite{STZ} for families of $K3$-surfaces, and methods and results
of \cite{STZ} have been our motivation to study the case of abelian varieties.

\begin{theorem}\label{geomsplitt}
Let $f:A\to Y$ be a family of abelian varieties smooth over $Y\setminus S$,
and such that the local monodromies around $s\in S$ are unipotent.
If $S\not=\emptyset,$ and if $f:A\to Y$
reaches the Arakelov bound, then there exists an \'etale covering
$\pi:Y'\to Y$ such that $f':A'=A\times_YY'\to Y'$ is isogenous over $Y'$
to a product
$$ B\times E\times_{Y'}...\times_{Y'} E ,$$
where $B$ is abelian variety defined over $\C$ of dimension
$g-g_0,$ and where $h: E\to Y'$ is a family of elliptic
curves reaching the Arakelov bound.
\end{theorem}

We do not know whether for all $g$ there are families of Jacobians
among the families of abelian varieties considered in
\ref{geomsplitt}, i.e. whether one can find a family $\varphi:Z\to
Y$ of curves of genus $g$ such that $f:J(Z/Y)\to Y$ reaches the
Arakelov bound.

For $Y=\P^1$ the Arakelov inequality implies $\#S\geq 4$.
Our hope, that a family with $\# S=4$ can not be a family of
Jacobians, hence that the Jacobian of a family of curves over $\P^1$
must have more than $5$ singular fibres,
was destroyed by an example of a family of genus $2$ curves over the
modular curve $X(3)$ in \cite{Kan1}, whose Jacobian is isogenous to the
product of a fixed elliptic curve $B$ with the modular curve
$E(3)\to X(3)$ (see Section \ref{jacobians}).\\

As mentioned already, this note owes a lot to the recent work of
the second named author with Xiao-Tao Sun and Sheng-Li Tan.
We thank Ernst Kani for explaining his beautiful construction in \cite{Kan1},
and for sharing his view about higher genus analogs of families of curves with
splitting Jacobians. It is also a pleasure to thank H\'el\`ene Esnault for
her interest and help, and Ngaiming Mok, for explaining us differential geometric
properties of base spaces of families.

This note grew out of discussions during a visit of the first named author
in Hong Kong. He would like to thank the members of the
Institute of Mathematical Science and the Department of Mathematics at the
Chinese University of Hong Kong for their hospitality.

\section{Splitting of $\C$-local systems}

We will frequently use C. Simpson's correspondence between
polystable Higgs bundles of degree zero and representations of the
fundamental group $\pi_1(U,s)$.

\begin{theorem}[C. Simpson \cite{Si2}]\label{simpson1}
There exists a natural equivalence between the category of direct
sums of stable filtered regular Higgs bundles of degree zero, and
of direct sums of stable filtered local systems of degree zero.
\end{theorem}
We will not recall the definition of a ``filtered regular'' Higgs bundle
(\cite{Si2}, p. 717), just remark that for a Higgs bundle corresponding
to a local system with unipotent monodromy around the points in
$S$, the filtration is trivial.

For example, \ref{simpson1} implies that the splitting of Higgs
bundles (0.1.1) corresponds to a decomposition over
$\mathbb C$
$$ (R^1f_*\Z_{A_0})\otimes \C=\V \oplus\U_1$$
where $\V$ corresponds to the Higgs bundle $(F=F^{1,0}\oplus
F^{0,1}, \theta)$ and $\U_1$ to $(N=N^{1,0}\oplus N^{0,1},\theta_N=0)$.
 Let $\Theta(N,h)$  denote the curvature of the Hodge metric $h$ on
$E^{1,0}\oplus E^{0,1}$ restricted to $N,$
then by \cite{G}, chapter II we have
$$
\Theta(N,h|_N)=-\theta_N\wedge\bar\theta_N-\bar\theta_N\wedge\theta_N=0.
$$
This means that $h|_N$ is a flat metric. Hence,  $\U_1$ is a
unitary local system.

The local system $\V$ on $Y\setminus S$ is a variation of Hodge
structures with unipotent local monodromies around $s\in S.$ Hence by
Deligne's theorem \cite{Del}, one obtains a decomposition
\begin{equation}\label{splitt2}
\V =\bigoplus_i \mathbb V_i,
\end{equation}
for irreducible local systems $\mathbb V_i $.

Restricting the Hodge filtration of $\V$ to $\V_i$, one obtains a
Hodge filtration on $\V_i$, which in general is not polarized, and
(\ref{splitt2}) is a decomposition of $\C$-variations of Hodge
structures. Taking the grading of the Hodge filtration, one
obtains a decomposition of the Higgs bundle $(F=F^{1,0}\oplus
F^{0,1}, \theta)$ as a direct sum of sub Higgs bundles, as stated
in \ref{simpson1}.\\

As a typical application of Simpson's correspondence one finds
\begin{proposition}\label{semistable}
Keeping the notations from the Introduction, assume that $f:A\to
Y$ reaches the Arakelov bound. Then
\begin{enumerate}
\item[i)] The sheaf $F^{1,0}$ is  poly-stable. Namely there is a
decomposition
$$ F^{1,0}\simeq \bigoplus_i \sA_i$$
with $\sA_i$ stable, and
$$\frac{\deg \sA_i}{\rank\sA_i}=\frac{\deg F^{1,0}}{\rank F^{1,0}}.$$
Moreover, $\theta_{1,0}:F^{1,0} \to F^{0,1}\otimes \Omega^1_Y(\log S)$
is an isomorphism.
\item[ii)] If $\deg\Omega^1_Y(\log S)$ is even, there exists a
decomposition
$$ \V\simeq \L\otimes \T,$$
such that $\T$ is a unitary local system,
and $\L$ is a rank-2 local system. For some invertible
sheaf $\sL$ the Higgs bundle corresponding to $\L$ is
$(\sL\oplus\sL^{-1},\tau)$, with $\tau|_{\sL^{-1}}=0$ and
$\tau|_{\sL}$ given by an isomorphism
$$
\tau^{1,0}: \sL\to \sL^{-1}\otimes\Omega^1_Y(\log S).
$$
\end{enumerate}
\end{proposition}

\begin{proof} i) \quad Let $\sA\subset F^{1,0}$ be a subsheaf, and let
$\sB\otimes\Omega_Y^1(\log S)$ be its image under $\theta_{1,0}$.
Then $\sA\oplus\sB$ is a Higgs subbundle of $F^{1,0}\oplus
F^{0,1}$, and applying \ref{simpson1} one finds
$\deg(\sA)+\deg(\sB)\leq 0$. Hence
\begin{multline*}
\deg(\sA) = \deg(\sB) + \rank(\sB)\cdot\deg(\Omega_Y^1(\log S))\\
\leq -\deg(\sA) + \rank(\sA)\cdot\deg(\Omega_Y^1(\log S)),
\end{multline*}
and (\ref{arakelov_ineq}) implies that
$$
\frac{\deg(\sA)}{\rank(\sA)} \leq \frac{1}{2}\deg(\Omega_Y^1(\log S))=
\frac{\deg(F^{1,0})}{g_0}.
$$
By \ref{simpson1} the Higgs bundle $(F^{1,0}\oplus
F^{0,1},\theta)$ splits as a direct sum of stable Higgs bundles of
degree zero. If
$$
\frac{\deg(\sA)}{\rank(\sA)} = \frac{\deg(F^{1,0})}{g_0},
$$
the degree of $\sA\oplus\sB$ is zero, ${\rm rank}(\sA)={\rm rank}(\sB)$, and
$(\sA\oplus\sB,\theta|_{\sA\oplus\sB})$ is a direct factor of
$(F^{1,0}\oplus F^{0,1},\theta)$. In particular, $\sA$ is a direct
factor of $F^{1,0}$.

For $\sA=F^{1,0}$ one finds $\theta_{1,0}$ to be injective. By 
(\ref{arakelov_ineq}) it must be an isomorphism.\\
\ \\
ii)\quad  Taking the determinant of
$$
\theta^{1,0}: F^{1,0}\> \simeq >> F^{0,1}\otimes\Omega^1_Y(\log
S),
$$
one obtains an isomorphism
$$
\det\theta^{1,0}: \det F^{1,0}\> \simeq >> \det
F^{0,1}\otimes{\Omega^1_Y(\log S)}^{\otimes g_0},
$$
Since $F^{1,0}\simeq F^{0,1\vee},$
$$
(\det F^{1,0})^{\otimes 2}\simeq {\Omega^1_Y(\log S)}^{\otimes
g_0}.
$$
By assumption there exists an invertible sheaf $\sL$ with
$\sL=\Omega^1_Y(\log S)^{1/2}.$ For some invertible
sheaf $\sN$ of order two in ${\rm Pic}(Y)$, 
one finds $\sN\otimes \det F^{1,0}=\sL^{\otimes g_0}$. By part i) the sheaf
$$\sT=F^{1,0}\otimes \sL^{-1}$$
is poly-stable of degree zero. \ref{simpson1} implies that the
Higgs bundle $(\sT,0)$ corresponds to a local system $\T$,
necessarily unitary.

Choose $\L$ to be the local system corresponding to the Higgs
bundle
$$(\sL\oplus\sL^{-1},\tau), \mbox{ \ \ with \ \ }
\tau^{1,0}: \sL \> \simeq >> \sL^{-1}\otimes \Omega^1_Y(\log S).
$$
The isomorphism
$$
\theta^{1,0}:\sT\otimes \sL=F^{1,0}\> \simeq >> F^{0,1}\otimes
\Omega^1_Y(\log S)\> \simeq >> F^{0,1}\otimes \sL^{\otimes 2}
$$
induces an isomorphism
$$
\phi: F^{0,1}\> \simeq >> \sT\otimes \sL^{-1},
$$
such that $\theta^{1,0}={\rm id}_\sT \otimes \tau^{1,0}$. Hence
the Higgs bundles $(F^{1,0}\oplus F^{0,1},\theta)$ and
$(\sT\otimes(\sL\oplus \sL^{-1}),{\rm id}_\sT\otimes\tau)$ are
isomorphic, and $\V\simeq \T\otimes \L$.
\end{proof}

\begin{remark}\label{even} \
\begin{enumerate}
\item[i)] If $\deg \Omega^1_Y(\log S)$ is odd, hence $S\neq \emptyset$, and
if the genus of $Y$ is not zero, one has to replace $Y$ by an
\'etale two to one cover, in order to be able to apply
\ref{semistable}, ii).
\item[ii)]
For $Y=\P^1$ and for families reaching the Arakelov bound,
$\#S$ is always even. This, together with 
the decomposition \ref{semistable}, ii),
for $\U=\C^{g_0}$, can easily obtained in
the following way. By \ref{semistable}, i), $F^{1,0}$ must be the
direct sum of invertible sheaves $\sL_i$, all of the same degree,
say $\nu$. Since $\theta^{1,0}$ is an isomorphism, the image
$\theta^{1,0}(\sL_i)$ is $\sO_{\P^1}(2-s+\nu)\otimes \Omega$.
Since $F^{0,1}$ is dual to $F^{1,0}$ one obtains $-\nu=2-s+\nu$,
and writing $\sL_i^{-1}=\theta^{1,0}(\sL_i)$,
$$
(F^{1,0}\oplus F^{0,1}, \theta)\simeq (\bigoplus_i
\sO_{\P^1}(\nu)\oplus \sO_{\P^1}(-\nu), \bigoplus_i \tau).
$$
\end{enumerate}
\end{remark}

Consider now the endomorphism $\sEnd(\V)$ of $\V$, which is a
weight zero variation of Hodge structures. The Higgs bundle
$$
(F^{1,0}\oplus F^{0,1},\theta)
$$
for $\V$ induces the Higgs bundle
$$
(F^{1,-1}\oplus F^{0,0}\oplus F^{-1,1},\theta)
$$
corresponding to $\sEnd(\V)=\V\otimes\V^\vee$, by choosing
\begin{gather*}
F^{1,-1}=F^{1,0}\otimes {F^{0,1}}^\vee,
\ \ \ F^{0,0}=
F^{1,0}\otimes {F^{1,0}}^\vee \oplus F^{0,1}\otimes {F^{0,1}}^\vee\\
\mbox{ \ \ \ and \ \ \ }
F^{-1,1}=F^{0,1}\otimes {F^{1,0}}^\vee.
\end{gather*}
The Higgs field is given by
$$
\theta_{1,-1}= -\id\otimes{\tau_{1,0}}^\vee \oplus \tau_{1,0}\otimes \id
\mbox{ \ \ \ and \ \ \ }
\theta_{0,0}=\tau_{1,0}\otimes \id \oplus -\id \otimes {\tau_{1,0}}^\vee.
$$
\begin{lemma}\label{splitt3}
Assume that $f:A\to Y$ reaches the Arakelov bound
or equivalently that the Higgs field of $\V$ is maximal. Let
$$
F^{0,0}_u := {\rm Ker} (\tau_{0,0}) \mbox{ \ \ \ and \ \ \ }
F^{0,0}_m = {\rm Im}(\tau_{1,-1}).
$$
Then there is a splitting of the Higgs bundle
$$
(F^{1,-1}\oplus F^{0,0}\oplus F^{-1,1},\theta)=
(F^{1,-1}\oplus F_m^{0,0} \oplus F^{-1,1},\theta)\oplus (F^{0,0}_u,0),
$$
which corresponds to a splitting of the local system over $\C$
$$
\sEnd(\V)=\W\oplus\U.
$$
Moreover, $\U$ is unitary of rank $g_0^2$, and $\W$ is a $\C$
variation of Hodge structures with maximal Higgs field, i.e.
$$
\tau_{1,-1}:F^{1,-1}\to F^{0,0}_m\otimes \Omega^1_Y(\log S)
\mbox{ \ \ \ and \ \ \ }
\tau_{0,0}:F^{0,0}_m\to F^{-1,1}\otimes \Omega^1_Y(\log S)
$$
are both  isomorphisms.
\end{lemma}
\begin{proof}
By definition, $(F^{0,0}_u,0)$ is a sub Higgs bundle of
$(F^{1,-1}\oplus F^{0,0}\oplus F^{-1,1},\theta)$. We have an exact
sequence
$$
0\to F^{0,0}_u \to F^{0,0} \to F^{-1,1}\otimes \Omega_Y^1(\log S) \to \sC
$$
where $\sC$ is a skyscraper sheaf. Hence
$$
\deg(F^{0,0}_u) \geq \deg(F^{0,0}) - \deg(F^{-1,1}) - \rank(F^{-1,1})\cdot
\deg(\Omega^1_Y(\log S).
$$
Note that if this inequality is an equality then $\sC$ is
necessarily zero. Since $\deg(F^{0,0})=0$ and since, by the
Arakelov equality,
$$
\deg(F^{-1,1}) = g_0\cdot \deg(F^{0,1}) + g_0 \cdot \deg({F^{1,0}}^\vee)
= g_0^2\cdot \deg(\Omega_Y^1(\log S))
$$
one finds $\deg(F^{0,0}_u) \geq 0$. By \ref{simpson1} the degree
of $F_u^{0,0}$ can not be strictly positive, hence it is zero.
Again by \ref{simpson1} $(F^{0,0}_u,0)$ being a Higgs subbundle of
degree zero with trivial Higgs field, it corresponds to a unitary
local subsystem $\U$ of ${\rm End}(\V)$. The exact sequence
$$
0\to F^{0,0}_u\to F^{0,0}\to F^{-1,1}\otimes\Omega^1_Y(\log S)\to 0
$$
splits, and one obtains a direct sum decomposition of Higgs
bundles
$$
(F^{1,-1}\oplus F^{0,0}\oplus F^{-1,1},\theta)= (F^{1,-1}\oplus
F_m^{0,0} \oplus F^{-1,1},\theta)\oplus (F^{0,0}_u,0),
$$
which induces the splitting on $\sEnd(\V)$ with the desired
properties.
\end{proof}

\begin{remark}\label{altsplitt}
Using \ref{semistable}, ii), the splitting in \ref{splitt3}
can be made more precise.
We know that $\V= \T\oplus\L$
with $\T$ unitary and $\L$ a weight two
variation of Hodge structures with maximal Higgs field.
One obtains
$$
\sEnd(\V)= \T\otimes {\T}^\vee \otimes \L \otimes {\L}^\vee.
$$
Applying \ref{splitt3} to $\L$ instead of $\V$, we obtain
a decomposition
$$
\sEnd(\L)=\L \otimes {\L}^\vee = \C \oplus {\mathbb S}
$$
where the $\C$ factor acts by multiplication on $\L$
and where ${\mathbb S}$ has a maximal Higgs field.
So $\T\otimes {\T}^\vee$ is a direct factor of $\V$
of rank $g_0^2$. Its complement
$\W=\T\otimes {\T}^\vee\otimes {\mathbb S}$ has again a maximal Higgs field.
\end{remark}

\begin{remark}\label{wedge}
If one replaces ${\rm End}(\V)$ by the isomorphic locally constant
system $(\V \otimes \V)\otimes_\Z\C$, one obtains the same
decomposition. However, it is more natural to shift the weights by
two, and to consider this as a variation of Hodge structures of
weight $2$.

A statement similar to \ref{splitt3} holds true for
$\wedge^2(\V)$. Here the Higgs bundle is given by
\begin{gather*}
{F'}^{2,0}=F^{1,0}\wedge F^{1,0},
\ \ \ F^{1,1}=
{F'}^{1,0}\otimes F^{0,1}
\mbox{ \ \ \ and \ \ \ }
{F'}^{0,2}=F^{0,1}\wedge F^{0,1}.
\end{gather*}
\end{remark}

\section{Splitting over $\bar{\Q}$}

Up to now, we tried to describe the local systems of $\C$-vector
spaces $\V$ induced by the family of abelian varieties. We say
that such a local system $\V$ is defined over a subfield $K$ of
$\C$, if there exists a local system $\V_K$ of $K$-vector spaces
with $\V=\V_K\otimes_K\C$. In this section we want to show, that
the splitting $\V=\W\oplus\U$ considered in the last section are
defined over $\bar{\Q}$, i.e. that there exists a number field $K$
and local systems $\V_K$, $\W_K$ and $\U_K$ with
$$
\V=\V_K\otimes \C, \ \ \W=\W_K\otimes \C, \ \ \U=\U_K\otimes
\C, \mbox{ \ and with \ } \V_K=\W_K \oplus \U_K.
$$
We start with a simple observation. Suppose that $\V$ is a local
system defined over a number field $K$. Choosing a base point $p\in
Y\setminus S$ the local system $\V_{K}$ is given by a
representation $\rho: \pi_1(Y\setminus S,p) \to
\text{Gl}(V_{K})$ for the fibre $V_{K}$ of $\V_{K}$ over $p$.

Fixing a positive integer $r,$ let $\mathcal G(r,\V)$ denote the
set of all rank-r sub local systems of $\V$. Then $\mathcal
G(r,\V)$ is the subvariety of the Grassmann variety ${\rm
Grass}(r, V_{K})$ consisting of the $\pi_1(Y\setminus S,p)$
invariant points. In particular, it is a projective variety
defined over $K$. An $L$-valued point of $\mathcal G(r,\V)$
corresponds to a sub local system of $\V_L=\V_{K}\otimes_K L$. One
obtains the following well known property.

\begin{lemma}\label{fieldofdef}
If $[\W]\in \mathcal G(r,\V)$ is an isolated point,
then $\W$ is defined over $\bar\Q.$
\end{lemma}

\begin{lemma}\label{splittqbar}
Let $\V$ be the underlying local system of an variation of Hodge
structures defined over a real number field $K$, and suppose that
there is a splitting
\begin{equation}\label{splitt5}
\V=\W\oplus\U,
\end{equation}
such that $\U$ is a unitary, and such that $\W$ has a generically
maximal Higgs field $\theta^{p,q},$ i.e
$$ \theta^{p,q}: F^{p,q}_{\W}\to F^{p-1,q+1}_{\W}\otimes\Omega^1_Y(\log S) $$
is generically isomorphic for all $(p,q)$ with $p>0$. Then this
splitting can be defined over $\bar\Q\cap \R,$  and is orthogonal
with respect to the polarization.
\end{lemma}

\begin{proof} Consider a family $\{\W_t\}_{t\in \Delta}$ of local subsystems of
$\V$ defined over a disk $\Delta$ with $\W_0=\W.$ For $t\in
\Delta$ let $(F_{\W_t},\theta_t)$ denote the Higgs bundle
corresponding to $\W_t$. Hence $(F_{\W_t},\theta_t)$ is obtained
by restricting the $F$-filtration of $\V\otimes \sO_U$ to
$\W_t\otimes\sO_U$ and by taking the corresponding graded sheaf.
So the Higgs map
$$ \theta^{p,q}: F^{p,q}_t\to F^{p-1,q+1}_t\otimes\Omega^1_Y(\log S)$$
will again generically isomorphic for $t$ sufficiently closed to
$0$ and $p>0$. If the projection
$$\rho: \W_t \to \V \to \W\oplus\U\to \U$$
is non-zero, the complete reducibility of local systems coming
from variations of Hodge structures (due to Deligne [1])
implies that $\W_t$ contains a non-trivial unitary direct factor,
say $\U_t$. Restricting again the $F$ filtration and passing to
the corresponding graded sheaf, we obtain a non-trivial splitting
sub Higgs bundle $(F_{\U_t},0)$ of $(F_{\W_t},\theta_t)$,
contradicting the generic maximality of the Higgs field $\oplus
\theta^{p,q}$. Hence $\rho$ is zero and $\W_t=\W$.

Thus $\W$ is rigid as a sub local system of $\V$, and by Lemma
\ref{fieldofdef} $\W$ is defined over $\bar{\Q}$.

By assumption $\V=\V_\R\otimes \C$ and the complex conjugation
defines an involution $\iota$ on $\V$. Let $\bar\W$ denote the
image of $\W$ under $\iota$. Then $\bar\W$ has again generically
isomorphic Higgs maps $\theta^{p,q},$ for $p>0$. If $\bar\W\neq\W$,
repeating the argument used above, one obtains a non-trivial map
$\bar \W\to \U$, contradicting again the generic maximality of the
Higgs field.

So enlarging the real algebraic number field $K$, if necessary, we
may assume that $\W=\W_K\otimes_K \C$ for some $\W_K \subset
\V_K$. The polarization on $\V$ restricts to a non-degenerated
intersection form on $\V_K$. Choosing for $\U_K$ the orthogonal
complement of $\W_K$ in $\V_K$ we obtain a splitting
$$ \V_K=\W_K\oplus\U_K$$
inducing over $\C$ the one in (\ref{splitt5}).
\end{proof}

\begin{corollary}\label{splitt6}
Let $f:A\to Y$ be a family of abelian varieties reaching the
Arakelov bound, and let $\U_1$ be a unitary local system with
$$ R^1f_*(\Z_{A_0})\otimes\C=\V\oplus\U_1,$$
such that $\V$ has a maximal Higgs field.
\begin{enumerate}
\item[i)] Then this splitting can be defined over $\bar\Q \cap
\R,$ and is it orthogonal with respect to the polarization.
\item[ii)] The splitting $\sEnd(\V)=\W\oplus\U$ constructed in
Lemma \ref{splitt3} can be defined over $\bar\Q\cap\R,$ and is
orthogonal with respect to the polarization.
\end{enumerate}
\end{corollary}

\section{Splitting over $\Q$}

\begin{lemma}\label{splitt7}
Assume that $S\neq \emptyset$ and let $\V_\Q$ be a $\Q$-variation
of Hodge structures of weight $k$. Assume that over some number
field $K$ there exists a splitting
$$
\V_K=\V_\Q\otimes_\Q K=\W_K\oplus\U_K
$$
where $\U=\U_K\otimes_K\C$ is unitary and where the Higgs field of
$\W=\W_K\otimes_K\C$ is maximal. Then $\W$, $\U$ and the
decomposition $\V=\W\oplus\U$ are defined over $\Q$.
Moreover, $\U$ extends to a local system over $Y$.
\end{lemma}
\begin{proof}
Let $\T$ be a sub local constant system of $\W$. Writing
$$
\big( \bigoplus_{p+q=k}F_{\T}^{p,q}, \bigoplus_{p+q=k}
\theta_{p,q}\big),
$$
for the Higgs bundle corresponding to $\T$, the maximality of the
Higgs field for $\W$ implies that the Higgs field for $\T$ is
maximal, as well. In particular, for all $s\in S$ and for $p>0$
the residue maps
$$
{\rm res}_s(\theta_{p,q}) : F_{\T,s}^{p,q} \>>> F_{\T,s}^{p-1,q+1}
$$
are isomorphisms. By \cite{Si2} the residues of the Higgs field at
$s$ are defined by the nilpotent part of the local monodromy
matrix around $s$. Hence if $\gamma$ is a small loop around $s$ in
$Y$, and if $\rho_\T(\gamma)$ denotes the image of $\gamma$ under
a representation of the fundamental group, defining $\T$, the
nilpotent part $N(\rho_\T(\gamma))=\log\rho_\T(\gamma)$ of $\rho_\T(\gamma)$
has to be non-trivial

We may assume that $K$ is a Galois extension of $\Q$. The Galois
group ${\rm Gal}(K/Q)$ acts on $\V_K$. For $\sigma \in {\rm
Gal}(K/Q)$ consider the composite
$$
p:\sigma(\U_K) \to \V_K=\W_K\oplus \U_K \to \W_K,
$$
and the induced map $\sigma(\U)=\sigma(\U_K)\otimes_K\C \to \W$.

Let $\gamma$ be a small loop around $s\in S$, and let
$\rho_\U(\gamma)$ and $\rho_{\sigma(\U)}$ be the images of
$\gamma$ under the representations defining $\U$ and $\sigma(\U)$
respectively. Since $\U$ is unitary and unipotent, the nilpotent
part of the monodromy matrix $N(\rho_{\U}(\gamma))=0$. This being
invariant under conjugation, $N(\rho_{\sigma(\U)}(\gamma))$ is
zero, as well as $N(\rho_{p(\sigma(\U))}(\gamma))$.

Therefore $p(\sigma(\U))=0$, hence $\sigma(\U)=\U$, and $\U$ is
defined over $\Q$. Taking again the orthogonal complement, one
obtains the $\Q$-splitting asked for in \ref{splitt7}.

Since $N(\rho_{\U}(\gamma))=0$, the residues of $\U$ are zero in
all points $s\in S$, hence $\U$ extends to a local system on $Y$.
\end{proof}

\begin{corollary}\label{splitt8} Suppose that $S\neq\emptyset.$ Then
the splittings in Corollary \ref{splitt6} can be defined over
$\Q.$
\end{corollary}

\section{$\Z$-structures and isogenies}

\begin{proposition}\label{splitt9}
Let $f:A\to Y$ be a family of abelian varieties with unipotent local monodromies
around $s\in S$, and reaching the Arakelov bound. If $S\not=\emptyset$ there exists a
finite \'etale cover $\pi:Y'\to Y$  with
\begin{enumerate}
\item[i)]
$\pi^*(R^1f_*(\Z_{A_0}))\supset \V'_{\Z}\oplus \Z^{g-g_0},\quad
\pi^*(R^1f_*(\Z_{A_0}))\otimes\Q=(\V'_{\Z}\oplus \Z^{g-g_0})\otimes\Q,$\\[.1cm]
where $\V'_{\Z}$ is an $\Z$-variation of Hodge structures of
weight $1$ with maximal Higgs field.
\item[ii)]
$\sEnd(\V'_{\Z})\supset \W'_{\Z}\oplus\Z^{g_0^2},\quad
\sEnd(\V'_{\Z})\otimes\Q=(\W'_{\Z}\oplus\Z^{g_0^2})\otimes\Q,$
where $\W'_{\Z}$ is an $\Z$-variation of Hodge structures of weight 0
with maximal Higgs field, and where $\Z^{g_0^2}$ is a constant
$\Z$-sub local system of type $(0,0)$.
\end{enumerate}
\end{proposition}

\begin{proof}
i) By \ref{splitt8} we already have the $\Q$-splitting
$$ R^1f_*(\Z_{A_0})\otimes\Q=\V_\Q\oplus\U_{1\Q}.$$
A $\Z$-structure on $\V_\Q$ and $\U_{1\Q}$ can be defined by
$$\V_{\Z}=R^1f_*(\Z_{A_0})\cap\V_\Q, \quad \U_{1\Z}=R^1f_*(\Z_{A_0})
\cap\U_{1\Q}.$$
Obviously
$$\V_{\Z}\otimes \Q=\V_\Q, \mbox{ \ \ and \ \ } \U_{1\Z}\otimes\Q=
\U_{1\Q}.$$
Since $\U_1$ is unitary and admits a $\Z$-structure, the monodromy
group of $\U_1$ is a finite group. Since the local monodromies of $\U_1$
around $S$ are trivial, $\U_1$ extends to a local system
$$ \rho_{\U_1}: \pi_1(Y,p)\to \text{Gl}(U_1),$$
where $U_1$ is the fibre of $\U_1$ in $p$.
After passing to the finite
\'etale cover of $\pi:Y' \>>> Y$ corresponds to $\rho_{\U_1}$  we obtain
a trivial $\Z$-sub local system of $\pi^*(R^1f_*(\Z_{A_0}))$ of
rank $g-g_0$. Together with $\pi^*\V_{\Z}$ we have
$$\pi^*(R^1f_*(\Z_{A_0}))\supset \pi^*\V_{\Z}\oplus \Z^{g-g_0},$$
such that
$$ \pi^*(R^1f_*(\Z_{A_0}))\otimes\Q=(\pi^*\V_{\Z}\oplus \Z^{g-g_0})
\otimes\Q.$$

ii) By and Lemma \ref{splitt6}, ii) and by Lemma \ref{splitt8} one
has a $\Q$-splitting
$$\sEnd(\V_{\Z})\otimes\Q=\W_\Q\oplus\U_\Q,$$
where $\U$ is a rank-$g_0^2$ unitary local system of (0,0)-type,
and where $\W$ has a maximal Higgs field. So ii) follows from the
same argument used to prove i).
\end{proof}

\begin{proof}[Proof of Theorem \ref{geomsplitt}]
Let $Y'$ be the \'etale covering constructed in \ref{splitt9},
ii). So using the notations introduced there,
\begin{equation}\label{zstructure}
R^1f'_*(\Z_{A'_0})\otimes \Q = V'_\Q \oplus \Z^{g-g_0}
\mbox{\ \ and \ \ }
\sEnd(\V'_\Q) = \W'_\Q \oplus \Z^{g_0^2}.
\end{equation}
The left hand side of (\ref{zstructure}) implies that $f':A'\to
Y'$ is isogenous to a product of a family of $g_0$ dimensional
abelian varieties with a constant abelian variety $B$. By abuse of
notations we will assume from now on, that $B$ is trivial, hence
$g=g_0$ and $R^1f'_*(\Z_{A'_0})\otimes \Q = V'_\Q$, and we will
show that under this assumption $f':A'\to Y'$ is isogenous to a
$g$-fold product of a modular family of elliptic curves.

Let us write $\End(*)=H^0(Y',\sEnd(*))$ for the global
endomorphisms. As explained in \cite{Sai}, for example,
$\End(\V'_\Q)=\Q^{g^2}$ is a $\Q$ Hodge structure of weight
zero, in our case the Hodge filtration is trivial, i.e.
$\End(\V'_\Q)^{0,0}=\End(\V'_\Q)$.

If $A_\eta=A'\times_{Y'}{\rm Spec}(\C(Y'))$ denotes the general
fibre of $f'$, one obtains
$$
\End(A_\eta)\otimes \Q = \End(\V'_\Q)^{0,0}=\End(\V'_\Q).
$$

By the complete reducibility of abelian varieties
there exists simple abelian varieties $B_1,\ldots,B_r$
of dimension $g_i$, respectively,
which are pairwise non isogenous, and such that
$\A_\eta$ is isogenous to the product
$$
B_1^{\times \nu_1} \times \cdots \times B_r^{\times \nu_r}.
$$
Moreover, since $\V$ has no flat part, none of the $B_i$
can be defined over $\C$. Let us assume that
$g_i=1$ for $i=1,\ldots ,r'$ and $g_i>1$ for
$i=r'+1,\ldots ,r$.

Let us write $D_i=\End(B_i)\otimes\Q$. By \cite{Mum}, p. 201,
Each $D_i$ is a division algebra of finite rank over $\Q$
with center $K_i$. Let us write $d_i^2$ for the rank of $D_i$
over $K_i$ and $e_i$ for the rank of $K_i$ over $\Q$.
Hence $e_i\cdot d_i^2$ is the rank of $D_i$ over $\Q$.

By \cite{Mum}, p. 202, or by \cite{BL}, p. 141, either $d_i \leq
2$ and $e_i \cdot d_i$ divides $g_i$, or else $e_i\cdot d_i^2$
divides $2\cdot g_i$. In both cases, the rank $e_i\cdot d_i^2$ is
smaller than or equal to $2\cdot g_i$. If $i\leq r'$, hence if
$B_i$ is an elliptic curve, not defined over $\C$, we have
$e_i=d_i=1$.

Writing $M_{\nu_i}(D_i)$ for the $\nu_i\times \nu_i$ matrices
over $D_i$, one finds (\cite{Mum}, p. 174)
$$
\End(A_\eta)\otimes \Q =
M_{\nu_1}(D_1)\oplus\cdots \oplus M_{\nu_r}(D_r)
$$
hence
\begin{multline*}
g^2= \dim_\Q(\End(A_\eta)\otimes \Q)=
\big(\sum_{i=1}^r \nu_i\cdot g_i \big)^2 =
\sum_{i=1}^r (e_i\cdot d_i^2) \cdot \nu_i^2
\leq\\
\sum_{i=1}^{r'}\nu_i^2 + \sum_{i=r'+1}^r \nu_i^2\cdot 2\cdot g_i
\leq \sum_{i=1}^r \nu_i^2\cdot g_i^2 .
\end{multline*}
Obviously this implies that $r=1$ and that $g_1 \leq 2$. If
$g_1=1$, we are done. In fact, the isogeny extends all over
$Y\setminus S$ and, since we assumed the monodromies to be
unipotent, $E:=B_1$ is the general fibre of a semi-stable family of
elliptic curves. The Higgs field for this family is
again maximal.

Before excluding the case $g_1=2$, let us compare this construction with
the one in Remark \ref{even}, for $Y=\P^1$ and in
Remark \ref{altsplitt} in general.

\begin{remark}\label{altsplitt2}
Writing $\T'$, $\L'$ and ${\mathbb S}'$ for the pullbacks of
$\T$, $\L$ and ${\mathbb S}$, respectively, one finds as remarked
in \ref{altsplitt} decompositions
$$
\sEnd(\V')= \T'\otimes {\T'}^\vee \otimes \L' \otimes {\L'}^\vee.
$$
and
$$
\sEnd(\L')=\L' \otimes {\L'}^\vee = \C \oplus {\mathbb S}'.
$$
So $\T'\otimes {\T'}^\vee$ is a direct factor of $\V$
of rank $g_0^2$, by construction unitary. Its complement
$\T'\otimes {\T'}^\vee\otimes {\mathbb S}'$ has again a maximal
Higgs field, hence no global section, and
$\T'\otimes {\T'}^\vee=\C^{g_0^2}$.
Altogether one finds
$$
\sEnd(\V')=\sEnd(\T')=\bigoplus^{g_0^2} \sEnd(\L').
$$
Using this description, one finds again that $\A_\eta$ can not be
the product of different non-isogenous abelian varieties.
\end{remark}
\noindent {\it End of the proof of \ref{geomsplitt}.} It remains to
exclude the case that $g_1=2$, and that $e_1\cdot d_1^2=4$. If the
center $K_1$ is not a totally real number field, $e_1$ must be
lager than $1$ and one finds\\[.1cm]
I. $d_1=1$ and $D_1=K_1$ is a quadratic imaginary extension of a
real quadratic extension of $\Q$.\\[.1cm]
If $K_1$ is a real number field, looking again to the
classification of endomorphisms of simple abelian varieties in
\cite{Mum} or \cite{BL}, one finds that $e_1$ divides $g_1$, hence
the only possible case is \\[.1cm]
II. $d_1=2$ and $e_1=1$, and $D_1$ is a quaternion algebra over
$\Q$.\\[.1cm]
The abelian surface $B_1$ over ${\rm Spec}(\C(Y'))$ extends to a
non-isotrivial family of abelian varieties $B'\to Y'$, smooth outside of
$S$ and with unipotent monodromies for all $s\in S$. This family
again has a maximal Higgs field, and thereby the local monodromies
in $s \in S$ are non-trivial. As we will see below, in both cases,
I and II, the moduli scheme of abelian surfaces with the corresponding type of
endomorphisms turns out to be a compact subvariety of the
moduli scheme of polarized abelian varieties, a contradiction.\\[.1cm]
I. By \cite{BL}, Example 6.6 in Chapter 9, there are only finitely
many $g_1$ dimensional abelian varieties with a given type of
complex multiplication, i.e. with $D_1$ a quadratic imaginary
extension of a real number field of degree $g_1$ over $\Q$.\\[.1cm]
II. By \cite{BL}, Exercise (1) in Chapter 9, there is no abelian
surface for which $D_1$ is a totally definite quaternion algebra.
Hence it remains to show the compactness of the moduli scheme of
abelian surfaces $B$ with a totally indefinite $D_1=\End(B)\otimes
\Q$, i.e. of the moduli scheme of false elliptic curves. Such
abelian surfaces and there moduli have been studied in \cite{Sha},
and there it is shown, that the moduli scheme is a compact Shimura
curve. This also follows from the construction of the moduli
scheme in \cite{BL}, \S 8 in Chapter 9, as a quotient of the upper
half plane $\sH$, and from \cite{Shi}, Chapter 9.
\end{proof}
\begin{proof}[Proof of Proposition \ref{modular}]
Let $\pi:Y'\to Y$ be an \'etale covering, $S'=\pi^{-1}(S)$ and let
$g:E\to Y'$ be a semi-stable family of elliptic curves, reaching
the Arakelov bound, and with $E_0=g^{-1}(Y'\setminus S')$ smooth,
for example the family occurring in \ref{geomsplitt}. Hence
$\L_{\mathbb Z}=R^1g_*\mathbb Z_{E_0}$ is a $\Z$-variation of Hodge
structures of weight one and of rank two.
Writing $\sL$ for the $(1,0)$ part, we have an isomorphism
\begin{equation} \label{modular1}
\tau_{1,0}: \sL \>>> \sL^{-1} \otimes \Omega_{Y'}(\log S').
\end{equation}
Since $\sL$ is ample, $\Omega_{Y'}(\log S')$ is ample, hence the
universal covering of $U'=Y'\setminus S'$ is the upper half plane
$\sH$. One obtains a commutative diagram
$$
\begin{CD}
\sH \> \tilde \varphi >> \sH\\
\V \psi' VV \V \psi VV\\
U' \> j >> \C
\end{CD}
$$
where $j$ is given by the $j$-invariant of the fibres of $E_0 \to
U'$, where $\psi$ is the quotient map $\sH \to
\sH/\text{Sl}_2(\Z)$, and where $\tilde\varphi$ is the period map.
Since the tangent sheaf of the period domain $\sH$ is given by the
sheaf of homomorphisms from the $(1,0)$ part to the $(0,1)$ part
of the variation of Hodge structures, the isomorphism $\tau_{1,0}$
implies that $\tilde \varphi$ is a local diffeomorphism. Note that
the Hodge metric on the Higgs bundle corresponding  to $\L_{\Z}$ has
logarithmic growth at $S$ and bounded curvature by Schmid \cite{Sch}.
By the remark after Prop. 9.8 together with the remark after Prop. 9.1
in \cite{Si3}
$$\tilde \varphi: \tilde U' {\to}\sH  $$
is a covering map, hence an isomorphism.

Since $\tilde \varphi$ is an equivariant isomorphism with respect
to the $\pi_1(U',*)-$action on $\tilde U'$ and the
$P\rho_{\L_\Z}(\pi_1(U',*))-$action on $\sH,$ the
homomorphism
$$\rho_{\L_\Z}: \pi_1(U',*)\to P\rho_{\L_\Z}(\pi_1(U',*))\subset \text{PSl}_2(\Z)$$
must be injective, hence an isomorphism.

This in turn implies, that
$$\varphi:  U'\to \sH/\rho_{\L_\Z}(\pi_1(U',*))     $$
is an isomorphism, $ \rho_{\L_\Z}(\pi_1(U',*))\subset \text{Sl}_2(\Z)$
is of finite index, and  $E_0\to U'$ is a modular curve.
\end{proof}

\section{Family of curves and Jacobians}\label{jacobians}

Let $Y$ be a curve, let $h: Z\to Y$ be a semi-stable
non-isotrivial family of curves of genus $g>1$, smooth over $V$,
and let $f:J(Z/Y)\to Y$ be a compactification of the Neron model of
the Jacobian of $h^{-1}(V)\to V$. Let us write $S$ for the points in $Y-V$
with $f^{-1}(y)$ singular and $\Gamma$ for the other points in
$Y\setminus V$, i.e. for the points $y$ with $h^{-1}(y)$ singular but
$f^{-1}(y)$ smooth. As usual we write $U=Y\setminus S$.

Let us first consider families of curves over $\P^1$.
S.-L. Tan \cite{Tan} has shown that $h:Z\to \P^1$
must have at least $5$ singular fibres, hence
$$
\# S + \# \Gamma
\geq 5.
$$
Moreover, he and Beauville \cite{Bea1} gave examples of
families with exactly $5$ singular fibres for all $g>1$. In those
examples one has $\Gamma = \emptyset$.

On the other hand, for $A=J(Z/Y)$ and for the ample sheaf
$F^{1,0}$ introduced in (\ref{splitt1}) the Arakelov inequality
(\ref{arakelov_ineq}) implies that
$$
2\cdot g_0 \leq 2\cdot\deg(F^{1,0}) \leq g_0 \cdot (\# S -2),
$$
hence $\# S \geq 4$. For $\# S = 4$, the family $f:J(Z/Y) \to Y$
reaches the Arakelov bound, hence by \ref{geomsplitt} it is
isogenous to a product of a constant abelian variety with a product
of modular elliptic curves, again with $4$ singular fibres.
By \cite{Bea2} there are just $8$ types of such families,
among them the universal family $E(3)\to X(3)$ of elliptic curves
with a level $3$-structures.

Being optimistic one could hope, that those families can not occur
as families of Jacobians, hence that there is no
family of curves $h:Z\to P^1$ with $\# S = 4$. However,
a counterexample has been constructed by E. Kani in \cite{Kan1}.

\begin{example}\label{elldiff}
Let $B$ be a fixed elliptic curve, defined over $\C$. Consider the
Hurwitz functor $\sH_{B,N}$ defined in \cite{Kan1}, i.e. the
functor from the category of complex schemes to the category of
sets with
\begin{multline*}
\sH_{B,N}(T)=\{f: C\to B\times T; \ f \mbox{ is a normalized
covering of degree }N \\
\mbox{ and }C\mbox{ a smooth family of curves of genus } 2 \mbox{
over } T \}.
\end{multline*}
The main result of \cite{Kan1} says that for $N \geq 3$ this functor is
represented by an open subscheme $V=H_{B,N}$ of the modular curve
$X(N)$ parameterizing elliptic curves with a level $N$-structure.

The smooth universal curve $\sC \to H_{B,N}$ extends to a semi-stable
curve $Z\to X(N)$ whose Jacobian is isogenous to $B\times E(N)$.
Hence writing $S$ for the cusps, $J(Z/X(N))$ is smooth outside
of $S$, whereas $Z\to X(N)$ has singular semi-stable fibres
outside of $H_{B,N}$. Theorem 6.2 in \cite{Kan1}
gives an explicit formula for the number of points in
$\Gamma=X(N)\setminus (H_{B,N}\cup S)$.

Evaluating this formula for $N=3$ one finds $\# \Gamma = 3$.
For $N=3$ the modular curve $X(3)$ is isomorphic to $\P^1$ with
$4$ cusps. So the number of singular fibres is $4$ for $J(Z/\P^1) \to \P^1$
and $7$ for $Z \to \P^1$.
\end{example}

We do not know whether similar examples exist for $g>2$. For $g>7$ the constant
part $B$ in Theorem \ref{geomsplitt} can not be of codimension one. In fact,
the irregularity $q(Z)$ of the total space of a family of curves of
genus $g$ over a curve of genus $q$ satisfies by
\cite{Xia}, p. 461, the inequality
$$
q(Z) \leq \frac{5\cdot g +1}{6} + g(Y).
$$
If $J(Z/Y) \to Y$ reaches the Arakelov bound, hence if
it is isogenous to a product
$$ B\times E\times_{Y}...\times_{Y} E ,$$
one finds $$\dim(B)\leq\frac{5\cdot g +1}{6}.$$

As explained in \cite{ES} it is not known, whether
for $g\gg 2$ there are any curves $C$ over $\C$ whose Jacobian
is isogenous to the product of elliptic curves.
We are even asking for families of curves
whose Jacobian is isogenous to the product of the same
elliptic curve, up to a constant factor.
\bibliographystyle{plain}

\end{document}